\documentclass[10pt]{amsart}
\usepackage[]{amsmath, amsthm, amsfonts, amssymb}
\usepackage{graphicx}
\usepackage[all]{xy}

\newcommand {\C} {{\mathbb C}}
\newcommand {\R} {{\mathbb R}}
\newcommand {\Z} {{\mathbb Z}}
\newcommand {\Q} {{\mathbb Q}}

\newcommand {\dt} {\bullet}
\newcommand {\F} {{\mathcal F}}

\newcommand {\cH} {{\mathcal H}}
\newcommand {\G} {{\mathcal G}}

\newcommand {\M} {{\mathcal M}}

\newcommand {\K} {{K\"ahler}}

\newtheorem{thm}[subsection]{Theorem}
\newtheorem{cor}[subsection]{Corollary}
\newtheorem{lemma}[subsection]{Lemma}
\newtheorem{prop}[subsection]{Proposition}

\newtheorem{ex}[subsection]{Example}

\begin{document}

\title{Homomorphisms between K\"ahler groups}

\author{ Donu Arapura}
 \address{
Department of Mathematics\\
 Purdue University\\
 West Lafayette, IN 47907\\
 U.S.A.}
 \thanks{Author partially supported by the NSF}

\maketitle

\begin{center}
  {\em To Anatoly Libgober}
\end{center}
This paper is an expanded version of  my talk at the Jaca conference; 
as such it is somewhere in between a  survey 
 and a research article. 
Algebraic geometry and topology have, of course, been connected from almost the beginning of both
subjects. The definition of the fundamental group is one of the first things one learns in a topology class, but
ironically it  appears to be one of the most subtle   and mysterious invariants of an  algebraic variety.
My talk was about the groups that arise as fundamental groups of projective manifolds or more generally
compact \K{} manifolds -- the so called \K{} groups.
The study of these groups came of age  in the 1980's and  90s.
During this period, tools ranging  from mixed Hodge theory, harmonic maps, rational homotopy theory and so on,
were applied to obtain strong restrictions on the groups in this class. In spite of all of the progress, some
difficult open questions remain. I will touch on a few  of these later on.

In this paper, I  want to propose one possible 
way forward by considering geometrically meaningful homomorphisms between K\"ahler  groups.
I define a homomorphism between such  groups to be K\"ahler (or 
K\"ahler-surjective)  if it comes from a holomorphic map (or a surjective holomorphic
map with connected fibres) between K\"ahler manifolds. 
Note that K\"ahler-surjective maps are surjective, but not all surjective K\"ahler homomorphisms are
 K\"ahler-surjective.
Some of the standard obstructions for a group to
be K\"ahler extend to K\"ahler homomorphisms in a fairly straightforward fashion.
A few examples are discussed below. However, I want to concentrate on 
a different sort of obstruction. I will say that a surjective homomorphism of groups $h:H\to G$ splits
if it has right inverse. There is a natural obstruction to splitting, namely the class $e(h)\in H^2(G, K/[K,K])$ associated
to the extension
$$0\to K/[K,K]\to H/[K,K] \to G \to 1 $$
where $K=\ker(h)$.  The main result here is that $e(h)$ vanishes in rationalized cohomology, when $h$
is K\"ahler-surjective. 
The vanishing result leads to a new  restriction for K\"ahler groups:
If a K\"ahler group
 admits a map $f$ onto a genus $g$ surface group $\Gamma_g$, with $g$ maximal, then $e(f)=0$ (rationally).
Thus for example, the group
$$\langle a_1,\ldots a_{2g},c\mid [a_1,a_{g+1}]\ldots
[a_{g},a_{2g}]=c, [a_i,c]=1\rangle
$$
is not  K\"ahler.

When the K\"ahler-surjective homomorphism comes from a surjective map $f:X\to Y$ of projective varieties,
there is a fairly simple  proof of the main theorem. This hinges on the fact
that $f$ has a multisection because the generic fibre has a rational point over a finite extension of
$\C(Y)$.  Of course this argument  is no longer valid in
the general K\"ahler setting,  so I use a different strategy. The first step is to show that $f:X\to Y$ can
be replaced by a map where $\ker\pi_1(f)$ has a reasonable structure.
There is a small price to be paid, in that $Y$ is now an orbifold.
Then the main step is to identify $e(\pi_1(f))$ as lying in the
image of the  differential of the Leray spectral sequence, and
then deduce the result from the  decomposition theorem
\cite{bbd, saito} suitably extended to K\"ahler orbifolds.

My thanks to J. Amor\'os, A. Suciu, and B. Wang for helpful comments.

\section{ Structure of K\"ahler groups}

I want to say a few words about what is and is not known  about the structure of \K{} groups,
while ignoring many other interesting  things pertaining to fundamental
groups of varieties (see \cite{hain}, \cite{libgober}, \cite{kollar}, \cite{schneps}, \cite{simpson2} and references therein).
My summary will be fairly sketchy.
The five author book \cite{abc} gives a much more complete account of the structure theory.

\subsection{What we know.}

 \begin{enumerate}
   \item (Hodge) The abelianization $\Gamma/D\Gamma,$ ($D\Gamma=[\Gamma,\Gamma]$)
 of a \K{} group has  even rank
by the Hodge decomposition. While elementary, this remark eliminates $50\%$ of all groups from the outset!
However, it doesn't say anything about a group such as the  Heisenberg group, which is a nontrivial 
 central extension of $\Z^2$ by $\Z$.

 \item (Deligne-Griffiths-Morgan-Sullivan \cite{dgms}) The rational homotopy type of a compact \K{} manifold is
as simple as possible. More precisely, it is ($1$-)formal. This implies that secondary operations such as
Massey products vanish. Furthermore  the Malcev Lie algebra of a  \K{} group 
$$L(\Gamma) = \frac{\Gamma}{[\Gamma,\Gamma]}\otimes\Q \times 
\frac{[\Gamma,\Gamma]}{[\Gamma,[\Gamma,\Gamma]]}\otimes \Q \times\ldots, $$
which is a kind of linearization, has a presentation with quadratic relations. This is quite a strong
restriction. One can see, for example, that the Malcev algebra of the Heisenberg group
is not quadratic, so it is not \K.
 \item (Goldman-Millson \cite{gm}) The singularities of the representation variety 
$Hom(\Gamma, GL_n(\C))$ of a \K{} group at points corresponding to semisimple
representations are  quadratic. This is closely related to the previous
result, and one can redo the Heisenberg group example from this point of view.

\item (A.-Nori \cite{an}, Delzant \cite{delzant})  Most solvable groups are not \K. The precise result of  Delzant
is that a solvable group is not \K{} unless it contains a nilpotent subgroup of finite index. The special case
for polycyclic groups, proved in \cite{an}, was based on a study of what has been variously called the
first cohomology support locus or jumping locus, or characteristic variety. I ought to mention that
this object has been studied quite extensively by Libgober and others.

\item (Carlson-Toledo \cite{ct1}, Simpson \cite{simpson})   
Recall that a discrete subgroup $\Gamma\subset G$ of a Lie group is a lattice if the Haar measure $vol(G/\Gamma)<\infty$.
While some lattices are  \K, many  are not.
These authors give specific obstructions for a lattice in 
a semisimple Lie group to be \K.
For example, lattices in $SO(1,n)$ or $SL(n,\R)$, with $n\ge 3$, are not \K.

\item (Gromov-Schoen \cite{gs}, A.-Bressler-Ramachandran \cite{abr}) ``Big'' groups, such as free products,
are not \K.
  \end{enumerate}

\subsection{What we don't know.}

Here are a few open problems in the area. Obviously, this list reflects my own taste and
biases.

\begin{enumerate}

\item (Johnson-Rees \cite{jr}) {\em Is the class of \K{} groups the same as the class of fundamental
groups of smooth projective varieties}? Inspired by work of Voisin \cite{voisin}, who constructed
compact \K{} manifolds which are not homotopic to projective manifolds, I will conjecture
the answer is no. Unfortunately, her examples don't shed any light here, since they do not have interesting fundamental
groups. In fact, I do not have any potential counterexamples;
  a possible candidate proposed in the introduction to \cite{an} has been taken out of the running by subsequent
work of Campana \cite{campana2} and Delzant \cite{delzant}.
I do want to mention as evidence that Botong Wang has recently found a counterexample which settles the analogous problem for \K{} homomorphisms. I will say more about this later.

\item It seems unreasonable to try  to characterize all \K{} groups; there is just
no way to get a handle on all of them.
A more reasonable goal would be try to {\em characterize the \K{} groups  within
a well understood subclass of all groups such as lattices in connected Lie groups}. This  can be specialized  further
to two subclasses with very distinct behaviours: Lattices in solvable groups and  lattices in semisimple groups.
\begin{enumerate}

\item By \cite{an} and some structure theory \cite{raghunathan},
 a lattice $\Gamma$ in a solvable group  is not \K{} unless the
ambient Lie group is nilpotent, i.e., unless $\Gamma$ is a torsion free finitely generated nilpotent
group. The restrictions on the  \K{} groups in this class coming from \cite{dgms} and other sources
are so severe  that only a few nonabelian examples are known \cite{campana,ct}. These are all
(larger rank) Heisenberg groups. 
It is not clear whether these examples are in any sense typical.
{\em So the problem is to either construct more examples of nilpotent \K{} groups or find more constraints
on them.}

\item {\em Carlson and Toledo \cite{ct1}  have conjectured that a   lattice in a semisimple Lie group $G$ is  not in
\K{} unless the associated symmetric space $G/K$ ($K=\text{a maximal compact}$) is Hermitian symmetric}. 
In more explicit terms, for the space to be Hermitian 
 the simple factors of $G$ must be one of 
$$SU(p,q), SO^*(2n), SO(p,2), Sp(n,\R)$$
or among a finite list of exceptional cases (see \cite{helgason}).
Note that the converse is often true.  That is if $G/K$  is a  Hermitian symmetric space, then 
a lattice  $\Gamma\subset G$ is \K{} if it is either  cocompact,  or arithmetic and the Baily-Borel
boundary of $\Gamma\backslash G/K$ has codimension at least $3$.
The results of \cite{ct1,simpson} mentioned above go part of the way toward this conjecture, 
but it appears that new ideas are needed. 

\end{enumerate}

\item {\em Is the genus $g$ mapping class group \K}?  I want to emphasize that this is open in spite of
 some misleading statements in the literature. The mapping class group can be viewed as the orbifold
fundamental group of the moduli space of curves $\M_g$.
An incorrect proof that this group  is \K{} is to take the Satake compactifcation
of $\M_g$, and cut by hyperplanes. But for this to work,
the codimension of the complement would have to be greater than $2$ (otherwise the fundamental group
changes when one cuts). However,
if one analyzes things carefully, one sees  that the Satake boundary does in fact have
a codimension $2$ stratum.  When $g=2$, Veliche \cite{veliche1,veliche2} showed that the mapping class group is not \K{} by reducing to the case of braid groups which were checked  to be non-\K{} in \cite{arapura}.  Veliche's argument doesn't generalize, so it is unclear what to expect for other genera. Farb \cite{farb}
suggests that these are \K{}, but I see the glass as half empty.

\item {\em Is the image of $H^*(\pi_1(X),\Z)\to H^*(X,\Z)$ a sub Hodge structure\footnote{I'm not 
sure who to attribute this to. It's come up at various times in conversations with Hain, Nori, and Srinivas.}}?
When $X$ is a variety, one can ask whether this is a submotive, and in particular
whether the image  $H^*(\pi_1(X),\Z_p)\to H^*(X,\Z_p)\cong H^*(X_{et},\Z_p)$ is invariant under the action
of $Aut(\C)$ on \'etale cohomology. 
The paper \cite{arapura2} was motivated by the Hodge structure question, but it doesn't directly address it.
To attempt to answer it one ask:
{\em Does there exist a good analytic/algebro-geometric model for the classifying
map $X\to B\pi_1(X)=K(\pi_1(X),1)$ when $X$ is a \K{} manifold/algebraic variety}?
I am purposely being vague about what a ``good model'' actually means; part of the
problem is to make this precise. In fact, a certain analytic model is used in this paper,
but it is not ``good'' as far as  the above question is concerned.
One thing   is clear that the $K(\pi_1(X),1)$ cannot simply be a manifold
or variety in general, because $\pi_1(X)$ might have infinite cohomlogical dimension.
As recent work of Dimca, Papadima and Suciu \cite{dps} shows, one cannot expect
such a simple model even if one allows
$\pi_1(X)$ to be replaced by a commensurable group.

\end{enumerate}

\section{Elementary properties of \K{} homomorphisms}

To reiterate a homomorphism of \K{} groups is  K\"ahler (or 
K\"ahler-surjective)  if it comes from a holomorphic map (or a surjective holomorphic
map with connected fibres) between K\"ahler manifolds. It seems an interesting question of
whether there are reasonable analogues of the statements in \S 1.1. For now, note that
any constraint for \K{} groups which comes from the existence of some functorial 
structure automatically generalizes to \K{} homomorphisms.
For example:

\begin{lemma}
  If $h:\Gamma_1\to \Gamma_2$ is K\"ahler then the image, kernel and cokernel
of the induced map $\Gamma_1/D\Gamma_1\to \Gamma_2/D\Gamma_2$ has
even rank.
\end{lemma}

\begin{proof}
If $h$ is realized by a holomorphic map $f:X_1\to X_2$ of compact \K{} manifolds,
the groups $\Gamma_i/D\Gamma_i=H_1(X_i)$ carry  Hodge structures of weight $-1$ which are preserved by $f_*$.
\end{proof}

The Malcev Lie algebra $L(\Gamma)$  
has a filtration by the lower central series
$$C^0L(\Gamma)=L(\Gamma),\quad C^{n+1}(\Gamma) = [L(\Gamma), C^nL(\Gamma)]$$
The quotients of $L(\Gamma)/C^n$ gives an inverse 
system of finite dimensional nilpotent Lie algebras, and $L(\Gamma)$ can be identified with the limit.

\begin{prop}
  If $h:\Gamma_1\to \Gamma_2$ is K\"ahler, then the induced map 
of Malcev Lie algebras $L(\Gamma_1)\to L(\Gamma_2)$ strictly preserves the lower central series.
That is $h(C^nL(\Gamma_1))= h(L(\Gamma_1))\cap C^nL(\Gamma_2)$.
\end{prop}

\begin{proof}
  By a theorem of  Morgan \cite{morgan},  $L(\Gamma_i)/C^N$ carries a functorial
system of mixed Hodge structures, such that the weight filtration is exactly the
lower central series. The proposition follows from the corresponding
strictness statement in mixed Hodge theory.
\end{proof}

A more explicit obstruction is given by:

\begin{cor}
  Suppose $h:\Gamma_1\to \Gamma_2$ 
is a homomorphism of \K{} groups   such $h(\Gamma_1)\subset D\Gamma_2$.
and that the induced map of Malcev Lie algebras is nonzero.
Then $h$ cannot be \K{}. 
\end{cor}

\begin{proof} This is a consequence of the fact that a
 map strictly preserving a filtration is zero if the induced
map on the associated graded is zero.
\end{proof}

I want to point out a serious limitation of the notion of K\"ahler morphism, which  is
 that the  class is not closed under composition.
Thus the collection of K\"ahler groups and morphisms do not form a subcategory of 
the category of groups.

\begin{ex}
The maps
$$\Z^2\stackrel{p}{\to}\Z^4\stackrel{q}{\to}\Z^2$$
given by the matrices
$$p= \begin{pmatrix}1 &0\\0&1\\ 0&0\\0&0\end{pmatrix},\quad
q=
\begin{pmatrix}
  1 &0&0&0\\ 0&0&1&0
\end{pmatrix}
$$
are both \K{}, but $q\circ p$ 
 is not since it has a rank one image.
\end{ex}

I want to end this section by explaining the analogue of problem (1) of \S1.2.
Call a group (respectively homomorphism) projective
if can be realized as the fundamental group (respectively homomorphism between fundamental
groups) of a smooth complex projective variety (respectively induced by a morphism of projective varieties).

\begin{thm}[B. Wang]
There exists a \K{} morphism between projective groups which is not \K.
  \end{thm}

The proof, which is based on Voisin's method \cite{voisin}, will be written up by him separately.

\section{Reduction to Orbifolds}

The rest of this paper will devoted to the proof of the  theorem stated in the introduction about the
vanishing of the splitting obstruction associated to a \K-surjective map $\pi_1(f)$.
One of the difficulties that needs to be dealt with is that the
kernel of $\pi_1(f)$ can be quite wild.
In this section, I want to show how to reduce it to something more manageable.

 Given a  surjective
homomorphism of groups $h:H\to G$, let $H^{(1)}(h)= H/D\ker(h)$ with canonical 
surjection $h^{(1)}:H^{(1)}(h)\to G$. For each $n>1$, define
$H^{(n)}(h)$ by the diagram
$$
\xymatrix{
 0\ar[r] & K/DK\ar[r]\ar[d]^{n\cdot} & H^{(1)}(h)\ar[r]^{h^{(1)}}\ar[d] & G\ar[r]\ar[d]^{=} & 1 \\ 
 0\ar[r] & K/DK\ar[r] & H^{(n)}(h)\ar[r]^{h^{(n)}} & G\ar[r] & 1
}
$$
where $K=\ker(h)$, $n\cdot$ is multiplication by $n$, and  the left hand square is a pushout. We will write $H^{(n)}$ if
$h$ is understood. Let $e(h)$ denote the class $H^1(G, K/DK)$ defined earlier.
Observe that $n\cdot e(h)=0$ if and only if $H^{(n)}\to G$ splits.

Call a continuous map connected if all of its fibres are connected.
 Fix a  proper connected holomorphic map  $f:X\to Y$ of connected 
complex manifolds. Then we have a surjection 
$\pi_1(f):\pi_1(X)\to \pi_1(Y)$. 

\begin{lemma}
  Given a commutative diagram of topological spaces
$$
\xymatrix{
 X'\ar[r]\ar[d]^{f'} & X\ar[d]^{f} \\ 
 Y'\ar[r]^{g} & Y
}
$$
where $f'$ is proper and connected. If $g$ induces an
isomorphism of fundamental groups, then $e(\pi_1(f'))=0$ (respectively, is
torsion) implies
$e(\pi_1(f))=0$ (respectively, is torsion). In particular, this holds if
$g$ is a bimeromorphic
map of complex manifolds, or
if $g$ is the inclusion of the complement of an analytic subset of codim
$\ge 2$ in a complex manifold.
\end{lemma}

\begin{proof}
  If  $ne(\pi_1(f'))=0$ then $\pi_1(X')^{(n)}\to
  \pi_1(Y')=\pi_1(Y)$ has a splitting. By composition, we get a splitting of
 $\pi_1(X)^{(n)} \to \pi_1(Y)$.
\end{proof}

A variant of the previous lemma will be needed later.

\begin{lemma}\label{lemma:genfinred}
  Given a commutative diagram of topological spaces 
$$
\xymatrix{
 X'\ar[r]\ar[d]^{f'} & X\ar[d]^{f} \\ 
 Y'\ar[r]^{g} & Y
}
$$
Such that 
$f,f'$ induce surjections on fundamental groups, and $g:Y'\to Y$ is a generically finite proper map
between oriented manifolds.
If $e(\pi_1(f'))$ is
torsion then $e(\pi_1(f))$ is  torsion. 
\end{lemma}

\begin{proof}
Suppose that  $\pi_1(X')^{(n)}\to   \pi_1(Y')$ has a splitting. Then
 $\pi_1(X\times_Y Y')^{(n)}\to   \pi_1(Y')$ has a splitting. This
implies that 
$$(\pi_1(X)\times_{\pi_1(Y)} \pi_1(Y'))^{(n)}\to \pi_1(Y)$$
splits.
Therefore the pullback  $g^*e(\pi_1(f))$ to $H^2(\pi_1(Y'), K/[K,K]\otimes \Q)$
is zero. Now for any rational local system $V$ on $Y$, the maps
labeled $\alpha$ and $\beta$ in the commutative diagram
$$
\xymatrix{
  H^2(\pi_1(Y), V)\ar[r]\ar[d]^{\alpha} & H^2(\pi_1(Y'), g^*V)\ar[d] \\ 
 H^2(Y, V)\ar[r]^{\beta} & H^2(Y',g^*V)
}
$$
are easily seen to be injective. For $\alpha$ this is a standard result. It  follows from the
Leray spectral sequence for the classifying map $k:Y\to B\pi_1(Y)$
and the fact that the fibres of $k$ are simply connected.
For $\beta$, a left inverse is given by $\frac{1}{\deg g} g_*$.
 Therefore $e(\pi_1(f))\otimes \Q$ vanishes.

\end{proof}

\begin{lemma}
Let $h:H\to F$ be a surjective group  homomorphism, and 
let $f:F\to G$ be another surjective homomorphism such that
$\ker(f)/D\ker(f)$ is a possibly infinite torsion group. If $e(h)$ is torsion,
then $e(f\circ h)$ is torsion.
\end{lemma}

\begin{proof}
 Let $K= \ker (f\circ h)$, $L= \ker(f)$
 and $P= F\times_G H$ with its canonical
 projections $p:P\to F$ and $q:P\to H$. Observe that there is a commutative diagram
$$
\xymatrix{
 H\ar[rd]\ar[rrd]^{id}\ar[rdd]_{h} &  &  \\ 
  & P\ar[d]^{p}\ar[r]^{q} & H\ar[d]^{f\circ h} \\ 
  & F\ar[r]^{f} & G
}
$$

  Suppose that $ne(h)=0$ then $h^{(n)}:H^{(n)}(h)\to F$ splits. Therefore $P^{(n)}\to
  F$ also splits since $h^{(n)}$ factors through it. Thus $ne(p)=0$. 
Note that $e(p)$ is the restriction of $e(f\circ
 h)$ to $H^2(F, K/DK)$.  Therefore this implies  the  image of
  $e(f\circ h)$  in $H^2(F, K/DK\otimes \Q)$ vanishes.

The Hochschild-Serre spectral sequence yields
an exact sequence
$$ H^0(G, H^1(L, K/DK\otimes \Q))\to 
H^2(G,K/DK\otimes \Q)\stackrel{r}{\to} H^2(F, K/DK\otimes \Q)$$
Since $L$ acts trivially on $K/DK$ and since $L/DL$ is torsion, 
$$H^1(L, K/DK\otimes \Q)=Hom(L/DL, K/DK\otimes \Q))= Hom(L/DL, K/DK\otimes \Q)=0$$
Therefore $r$ is injective. Thus $e(f\circ h)\otimes \Q=0$.

\end{proof}

Fix a reduced effective divisor $D=\cup D_i\subset Y$ with simple normal
crossings and positive integers $m_i$ along each component. In
\cite[thm 4.1]{mo}, it is shown how to construct a smooth Deligne-Mumford
stack with a morphism $p:Y^{orb}\to Y$ which is the minimal object for which $p^*D_i=m_iD_i$. We will need to understand the structure to some extent, and in particular to see why this extends to the analytic category.
So as to  avoid too many abstractions, we describe this using the
older orbifold language (c.f. \cite{cr,moerdijk}). Recall that an analytic orbifold is given by  a locally finite atlas $\{(U_i,G_i), \phi_{ij}\}$ consisting of open sets $U_i\subset
\C^n$ stable under finite groups $G_i$, and analytic gluing functions 
$$U_i\supseteq U_{ij}\stackrel{\phi_{ij}}{\longrightarrow} U_{ji}\subseteq U_j$$ 
respecting the actions of $G_{ij}=\{g\in G_i\mid g(U_{ij})\subseteq U_{ij}\}$.
This can be conveniently packaged by the groupoid 
$$\xymatrix{
\G_1=\coprod_{ij} G_{ij}\times U_{ij} \ar@<1ex>[r]^>>>>>{s}\ar@<-1ex>[r]_>>>>>{t}& \coprod_i U_i\ar[l]|>>>>> u}=\G_0$$
$$\G_1 \times_{s,\G_0,t} \G_1 \stackrel{\mu}{\to} \G_1\stackrel{\iota}{\to}\G_1 $$
with source  $s(g,y)_{ij} = (y)_i$, target $t(g,y)_{ij} = (g\phi_{ij}(y))_j$,
 unit $u(y)_i = (1,y)_{ii}$, and multiplication  and inversion 
induced by the corresponding  operations on the groups. The groupoid
defines the corresponding  analytic
Deligne-Mumford stack. 

To construct $Y^{orb}$,
choose an atlas $\{V_i,\phi_{ij}\}$ for the manifold $Y$, consisting of polydisks $|y_j|<\epsilon$,  such that
for each $i$ either $V_i\cap D=\emptyset$ or $D_j$ is defined
by $y_j=0$. In the first case, let $U_i=V_i$ with the old coordinates $z_j=y_j$, and $G_i=\{1\}$. In the second case, using new coordinates $z_j=y_j^{1/m_j}$,
set $U_i=\{|z_j|<\epsilon^{1/m_j}\}$ with 
$G_i = \prod
\Z/m_j\Z$ acting by $z_j\mapsto \exp(2\pi\sqrt{-1}/m_j)z_j$.
For $\phi_{ij}$ we use the old transition functions in the $y$'s rewritten in terms of the $z$'s.
This defines $Y^{orb}$ as an orbifold. For our purposes a holomorphic map of orbifolds is a map which can be described locally by maps $U_i\to U_i'$ and  homomorphisms $G_i\to G_i'$
together with the obvious compatibilities, or equivalently by a morphism
of analytic groupoids (such maps are called good or
strong in the above references). In particular, we
have a map $Y^{orb}\to Y$ defined by the map
of atlases  $U_i\to V_i,G_i\to 1$ given the change of variables.

Note that contrary
to first appearances, the topology of $Y$ and $Y^{orb}$ should be
regarded as different unless $m_i=1$. 
This begs the question of what we even mean by the topology of an orbifold.
The quickest answer is that we can build a topological space
$B\G$ by taking the geometric realization of the  nerve of
the above topological groupoid, which is the simplicial space
$$B\G_\dt= \G_1 \times_{s,\G_0, t} \G_1\times_{s,\G_0, t} \ldots \G_1$$
with face maps described in \cite{segal}.
This involves  choices, but the weak homotopy type of $B\G$ depends
only on $Y^{orb}$, and so we will denote this by $[Y^{orb}]$. When this construction is applied to a manifold, we recover its homotopy type. The fundamental group $\pi_1(Y^{orb}):= \pi_1([Y^{orb}])$ can be understood
in a more explicit fashion. It is defined
so that its (set, abelian group,...)-valued representations correspond to locally constant sheaves (of sets, abelian groups...)
on $Y^{orb}$, which are given by  collections of $G_i$-equivariant locally constant sheaves $L_i$ on
$U_i$  with gluing isomorphisms $\phi_{ij}^*L_i|_{U_i\cap U_j}\cong
L_j|_{U_i\cap U_j}$ subject to the cocycle condition. So for instance a disk $\Delta$ ``modulo'' $\Z/m\Z$ is not
contractible; it has fundamental group $\Z/m\Z$.
To work    out $\pi_1(Y^{orb})$ explicitly, choose simple
loops $\gamma_i$ around the components $D_i$ of $D$. Then
$\pi_1(Y)$ can be identified with $\pi_1(Y-D)$ modulo the group generated by the $\gamma_j$ and its conjugates. While
$\pi_1(Y^{orb})$ is the quotient of $\pi_1(Y-D)$
by the  subgroup generated by the conjugates of $\gamma_i^{m_i}$. 

\begin{lemma}
  The map $Y^{orb}\to Y$ induces a surjection from $\pi_1(Y^{orb})$ to $\pi_1(Y)$. The abelianization of the kernel of this map
is a (possibly infinite)  torsion group.
\end{lemma}

\begin{proof}
  The first statement is clear. The kernel is generated by conjugates of the
$\gamma_i$ which are torsion elements. Therefore the abelianized kernel  is torsion since it is generated by torsion elements.
\end{proof}

The following seems well known when the base is a curve, e.g. \cite[lemma 3]{cko}.

\begin{lemma}\label{lemma:cko}
Let $f:X\to Y$ be a surjective holomorphic map of complex manifolds.
  Suppose that the discriminant of $f$ is a smooth divisor $D$, and
that $f^{-1}D$ is a divisor with normal crossings  such that the
restriction of $f$ to  the intersections
of components are submersions over $D$. Let $m_i$ denote the greatest
common divisor of the multiplicities of the components of $f^{-1}D_i$,
and construct $Y^{orb}$ as above.
Let $y_0\in Y-D$. Then $\pi_1(f)$ factors  through a surjection $\phi:\pi_1(X)\to \pi_1(Y^{orb})$, and
 $\pi_1(f^{-1}(y_0))$ surjects onto $\ker\phi$.

\end{lemma}

\begin{proof}
The map $f$ is given locally by
\begin{equation}
  \label{eq:cat}
  y_i = x_1^{k_{i1}}\ldots x_n^{k_{in}}
\end{equation}
The map $f$ factors through a map $X\to Y^{orb}$ given in the above coordinates by
$$
  y_i^{1/m_i} = x_1^{k_{i1}/m_i}\ldots x_n^{k_{in}/m_i}
$$
where $m_i = gcd(k_{i1}, \ldots k_{in})$.
This induces a factorization $\pi_1(X)\to \pi_1(Y^{orb})\to \pi_1(Y)$.
Now consider the commutative diagram
$$
\xymatrix{& \ker(r_2)\ar[r]\ar[d]&\ker(r_3)\ar[d]&\\
 \pi_1(f^{-1}(y_0))\ar[d]^{r_1}\ar[r] & \pi_1(X-f^{-1}D)\ar[d]^{r_2}\ar[r]^{\psi} & \pi_1(Y-D)\ar[d]^{r_3}\ar[r] & 1 \\ 
 \ker(\phi)\ar[r] & \pi_1(X)\ar[r]^{\phi} & \pi_1(Y^{orb})  &
}
$$
The middle row is exact since $f$ is a fibration over $Y-D$.
Furthermore the arrows $r_2, r_3$  are surjective. Therefore $\phi$ is also
surjective. This also implies that any element of $\ker(\phi)$ can be lifted
to an element of $\psi^{-1}\ker(r_3)\subset \pi_1(X-f^{-1}D)$. Therefore to
prove surjectivity of $r_1$,  it suffices to prove that $\psi|_{\ker(r_2)}$
surjects onto $\ker(r_3)$. By definition $\ker(r_3)$ is generated by conjugates of
$\gamma_i^{m_i}$, where $\gamma_i$ is a loop around $y_i=0$, using the
coordinates of \eqref{eq:cat}. A simple loop $\delta_j$ around the divisor $x_j=0$ maps to $\gamma_i^{\pm k_{ij}}$ in $Y-D$. Since $gcd(k_{ij})=m_i$, we can lift $\gamma_i^{m_i}$
to  a word in the $\delta_j$'s. Therefore any conjugate of
$\gamma_i^{m_i}$ can be lifted to  an element of $\ker(r_2)$. 
  
\end{proof}

\section{Decomposition theorem for orbifolds}

In  \cite{bbd}, Beilinson, Bernstein, Deligne and Gabber  developed
the theory of perverse sheaves on algebraic varieties. These form an
abelian subcategory of the derived category, and basic examples
include the intersection cohomology complexes $IC(L)[\dim Z]$ associated
to locally constant sheaves $L$ defined on locally closed sets $Z\subseteq X$.
 (We find it convenient to index $IC(L)$ so that $IC(L)=L$ generically on $Z$.)
They  proved a basic result called the decomposition theorem:
If $L$ is a semisimple perverse sheaf of geometric origin then for any proper map
$$\R f_*L =  \bigoplus_j IC(M_j)[m_j]$$
 for some $m_j\in \Z$ and $M_j$. 
Saito \cite{saito} has shown that this holds when $f$ is a proper holomorphic
map of K\"ahler manifolds. This will play a crucial role
in the proof of our theorem, however we will need a slight extension to orbifolds.

Suppose that $G$ is a finite group acting on a complex manifold $U$.
A $G$-equivariant sheaf is a sheaf $\F$ equipped with a a collection of
isomorphisms $\phi_g:g^*\F\cong \F$, for each $g\in G$, such that $\phi_1=id$ and
\begin{equation}
  \label{eq:equivF}
  \xymatrix{
 h^*g^*\F\ar[r]^{h^*\phi_g}\ar[rd]_{\phi_{gh}} & h^*\F\ar[d]^{\phi_h} \\ 
  & \F
}\tag{*}
\end{equation}
commutes. For example, the pullback of a sheaf from $U/G$ is naturally $G$-equivariant. But, unless the action of $G$ is free, not every equivariant arises this way.
Let  $Sh_G(U)$ denote the  category of  $G$-equivariant
sheaves of $\Q$-vector spaces. Suppose that  $K$ is 
the stabilizer of
the $G$-action on $U$. This  will act on a $G$-equivariant $\F$ by 
sheaf automorphisms, and the invariants $\F^K$ is a naturally a $G/K$-equivariant
sheaf. This gives an exact functor $\Gamma_K:Sh_G(U)\to Sh_{G/K}(U)$.
If $H\subset G$ is subgroup, then we have an induction functor $ind_H^G:Sh_H(U)\to Sh_G(U)$ right adjoint to the obvious restriction. Given a homomorphism
$f:H\to G$, we let $ind_H^G$ denote the composition $ind_{f(H)}^G\circ\Gamma_{\ker(f)}$.
Suppose that $X$ is an orbifold given by an atlas $\{(U_i, G_i),\phi_{ij}\}$. By a  sheaf
on $X$, we mean a collection of $G_i$-equivariant sheaves $\F_i$ together with
isomorphisms $\phi_{ij}^*\F_i|_{U_i\cap U_j}\cong \F_j|_{U_i\cap U_j}$
subject to the cocycle condition. Let $D^+(X)$ denote the derived category of sheaves on $X$.
Given a map $f:X\to Y$ of orbifolds
the direct image $\R f_* :D^+(X)\to D^+(Y)$  can be constructed in
couple of equivalent ways.   Sheaves on $X$ correspond  to simplicial
sheaves on the nerve  $B\G_\dt$, with respect to a given atlas,
such that the structure maps are all isomorphisms.
The  direct image can then realized  as a direct image of simplicial
sheaves \cite{hodge3}.
Alternatively, here is an explicit recipe:
If $f$ is  given by  $U_i\to V_i, G_i\to H_i$, and sheaf $\F=\{\F_i\}$ is a sheaf on $X$.
  Replace $\F_i$ by its
Godement flasque resolution $\G(\F_i)$ which is equivariant,
then form $ind_{G_i}^{H_i}f_*\G(\F_i)$.
An object in the derived category of sheaves on $X$ is a perverse sheaf if it restricts to 
a perverse sheaf in the usual sense on each $U_i$. $IC$ also extends to this setting.

By a K\"ahler orbifold, we simply mean that an atlas
 $\{(U_i,G_i),\phi_{ij}\}$ can
be chosen so that each $U_i$ possess a K\"ahler form preserved by
the $\phi_{ij}$.
 We can assume, by averaging over the groups,
 that these forms are invariant. We therefore have a K\"ahler class in
$H^2(X,\R)$.
So now we can state:

\begin{thm}[Saito$+\epsilon$]\label{thm:decomp}
  Suppose that  $f:X\to Y$ is a proper  holomorphic map of orbifolds with $X$
K\"ahler.
Let  $L$  be a geometric perverse sheaf on $X$  which means that it is a 
direct summand
of $IC(R^i\pi_*\Q)[\dim X]$ for some proper surjective holomorphic map $\pi:X'\to X$  of  K\"ahler orbifolds which is proper
over an open set. Then 
$$\R f_*L =  \bigoplus_j IC(M_j)[m_j]$$
\end{thm}

Restricting $f$ to the complement  of the discriminant  $D\subset Y$
allows us to identify some of these $M_j$'s and deduce that

\begin{cor}\label{cor:decomp}
  The perverse Leray spectral sequence decomposes as
$$E_2^{ij}= H^i(Y, IC(R^jf_*L|_{Y-D})) \oplus (\text{irrelevant stuff})$$
and furthermore it degenerates at $E_2$.
\end{cor}

We recall the basic properties of Saito's theory of mixed
Hodge modules \cite{saito1, saito2, saito}.

\begin{enumerate}
\item To any complex manifold $X$, there is an abelian category
of mixed Hodge modules $MHM(X)$ with a full semisimple
 abelian subcategory $MH(X)$ of pure  Hodge modules. (We are including
polarizability in the definition of $MH$.)

\item There is a faithful exact forgetful map from $MHM(X)$ to the category of rational perverse sheaves.

\item Given a closed irreducible analytic subvariety $Z\subseteq X$,
a pure Hodge module is said to have strict support $Z$
if the support of  the perverse sheaves  corresponding to all of  its simple
factors are all exactly $Z$. Any Hodge module decomposes uniquely into a sum of
submodules with strict support along various subvarieties.

\item If $Z$ is as above, any polarizable variation of Hodge structure $L$ on a Zariski
open subset of $Z$
extends to an object of $MH(X)$ such that the underlying 
perverse sheaf is the intersection cohomology complex $IC(L)[\dim Z]$. Conversely,
any Hodge module with strict support $Z$ is of this form.

\item All standard sheaf theoretic operations extend to $MHM$ and are compatible
with the corresponding operations on perverse sheaves. 

\item If $f:X\to Y$ is proper holomorphic map of manifolds with $X$ K\"ahler,
then the perverse sheaves  ${}^p\cH^i\R f_*\Q$ lift to objects in
$MH(Y)$. Moreover the hard Lefschetz holds, i.e. cupping with the 
K\"ahler class induces 
$${}^p\cH^{-i}\R f_*\R\cong {}^p\cH^i\R f_*\R(i)$$
in $MH(Y)\otimes \R$; the twist $(i)$ can be ignored for our purposes.
\end{enumerate}

Given a manifold with a finite group action, an equivariant mixed Hodge module
can be defined to be a mixed Hodge module $\F$ with  isomorphisms $\phi_g:g^*\F\cong \F$
satisfying the above compatibilities \eqref{eq:equivF}. A mixed Hodge module on an orbifold
$X$  with atlas $\{(U_i, G_i),\phi_{ij}\}$, is given by  a collection of $G_i$-equivariant modules $\F_i$ together with isomorphisms $\phi_{ij}^*\F_i|_{U_i\cap U_j}\cong \F_j|_{U_i\cap U_j}$.
 It is fairly easy to see that $MHM(X)$
still forms an abelian category with a forgetful functor to the category
perverse sheaves, which is defined in a similar fashion.
 A subvariety $Z\subset X$ corresponds to compatible
family of $G_i$-invariant subvarieties $Z_i\subseteq U_i$.
A pure Hodge module has strict support along $Z$ if each component has
strict support along $Z_i$. We can see that this notion is well
defined and that strict support decompositions hold. Moreover
the perverse sheaves  corresponding to  Hodge modules with strict support along $Z$ are
of the form $IC(L)[\dim Z]$ for a generic local system $L$ on $Z$.
The remaining statements are also true but require a bit more  explanation.

\begin{prop}\label{prop:eqMH}
Let $G$ be a finite group.  If $f:X\to Y$ is a $G$-equivariant proper
holomorphic map of manifolds with $X$ K\"ahler. Then  ${}^p\cH^i\R f_*\Q$ lifts to 
an equivariant pure Hodge module. This Hodge module is compatible
 with respect to restriction along open immersions.
\end{prop}

\begin{proof}
   All the essential ideas are due to Bernstein and Lunts \cite{bl}.
We can find a connected algebraic variety $V$
with  free $G$-action, such that 
\begin{equation}
  \label{eq:acyclicV}
H^i(V,\Q) =0\text{ for }0<i\le  N\gg 0   
\end{equation}
(for example, by fixing an embedding $G\subset GL_{n+N}(\C)$ and taking the Stieffel
variety $V$ of $n$-frames in $\C^{n+N}$). There are two key points. First, since the action of $G$ on $Y\times V$ is
free, the category of equivariant  Hodge modules on $Y\times V$ can be identified with $MH((Y\times V)/G)$.
Secondly, under
the projection $p:Y\times V\to Y$, we have an embedding of
the category of equivariant (mixed) Hodge modules on $Y$ to that on $Y\times V$ given by $p^*$.  To see this, first note that $p_*p^*\F=\F$ and $R^ip_*p^*\F=0$ for
and any sheaf $\F$ and  $0<i\le N$ by \eqref{eq:acyclicV}. Therefore,
given a pair of perverse sheaves $L,M\in Perv(Y)\subset D^{[-\dim Y,0]}(Y)$,
$$Hom(L,M)\to Hom(p^*L,p^*M)$$
is injective because it has a left inverse  given by $\tau_{\le 0}\R p_*$.
Since the $Hom$'s for the category of equivariant Hodge modules
are contained in the $Hom$'s of the underlying perverse sheaves, the earlier
claim about embeddings follows. 

Now consider the diagram
$$
\xymatrix{
 X\ar[d]^{f} & X\times V\ar[l]_{p'}\ar[r]^{q'}\ar[d]^{f_V} & (X\times V)/G\ar[d]^{F} \\ 
 Y & Y\times V\ar[l]_{p}\ar[r]^{q} & (Y\times V)/G
}
$$
Then by \cite{saito} ${}^p\cH^i\R F_*\Q$ lifts naturally to $MH((Y\times V)/G)$.
By base change, we can identify
$$q^*({}^p\cH^i\R F_*\Q)= {}^p\cH^i\R f_{V*}\Q= p^*({}^p\cH^i\R f_*\Q)$$
So the module on the right inherits the structure of
an equivariant  Hodge module.

Given a commutative diagram
$$
\xymatrix{
 X'\ar[d]^{f'}\ar[r]^{H} & X\ar[d]^{f}\\
 Y'\ar[r]^{h} & Y
}
$$
of $G$-equivariant maps with $h,H$ open immersions,
the map $h^*{}^p\cH^i\R f_*\Q\to {}^p\cH^i\R f'_*\Q$ is clearly compatible
with Hodge module structures.
\end{proof}

\begin{cor}
  If $f:X\to Y$ is a proper holomorphic map of  orbifolds with $X$ K\"ahler,
 ${}^p\cH^i\R f_*\Q$ lifts to a Hodge module. 
\end{cor}

 Combining this with the above remarks.

\begin{cor}\label{cor:decomp}
  If $f:X\to Y$ is a proper holomorphic map of orbifolds with $X$ K\"ahler,
 ${}^p\cH^i\R f_*\Q$  can be expressed
as a sum $\oplus IC(M_j)[\dim Z_j]$
\end{cor}

\begin{lemma}\label{lemma:hardlef}
   $f:X\to Y$ is a proper holomorphic map of orbifolds with $X$ K\"ahler,
the hard Lefschetz theorem holds for ${}^p\cH^i\R f_*\R$.
\end{lemma}

\begin{proof}
  This follows immediately from Saito's result \cite{saito}.
\end{proof}

\begin{proof}[Proof of theorem \ref{thm:decomp}]
By extending scalars to $\R$ and applying the previous lemma together
with \cite{deligne}, we get
$$\R f_*\Q \cong \bigoplus {}^p\cH^i\R F_*\Q[-i]$$
Thus by corollary \ref{cor:decomp}, we
can rewrite $\R f_*\Q$ as a sum of intersection cohomology complexes up to shift.
This proves the result for $L=\Q$. The general case can be deduced by applying the
theorem to $f\circ\pi$ as in \cite{saito}.
\end{proof}
 
\section{Splitting theorem}

We are ready to state the main theorem.

\begin{thm}\label{thm:splitting}
  Suppose that $f:X\to Y$ is a proper connected holomorphic map of
complex manifolds, with $X$ K\"ahler.
Then the splitting obstruction $e(\pi_1(f))\in H^2(\pi_1(Y), K/DK)$
is torsion, where  $K=  \ker(\pi_1(f))$.
\end{thm}

\begin{cor}
  If $h$ is a K\"ahler-surjective homomorphism then $e(h)$ is torsion.
\end{cor}

We note that K\"ahler-surjective homomorphisms are indeed both K\"ahler and surjective,
but not conversely.   Furthermore, the hypothesis  of this corollary cannot be replaced by the weaker condition.

\begin{ex}
  Campana \cite{campana,ct} has shown that certain Heisenberg groups $\Gamma$, which
 are nontrivial extensions of $\Z^{2n}$ by $\Z$, are K\"ahler.
The natural projection $\alpha:\Gamma\to \Gamma/D\Gamma=\Z^{2n}$ is K\"ahler
but $e(\alpha)$ is not torsion.
\end{ex}

Before starting the proof, we
need to recall some standard facts about classifying spaces.
Given a discrete group $G$, we can identify $BG=K(G,1)$. If $X$
is a good topological space (e.g. a CW complex),  then there is
a canonical map $k: X\to B\pi_1(X)$, unique up to homotopy,
 classifying the universal cover
viewed as a principle bundle over $X$. Given an orbifold $X$,
we thus get a classifying map 
$k:[X]\to B\pi_1([X])=B\pi_1(X)$, where $[X]$ is the associated homotopy type. We can realize this as
in a more explicit fashion by a simple modification
of the procedure given in \cite[\S 4]{segal}. 
Choose an atlas $\{U_i\}$ so that the intersections  $U_I=\cap_{i\in I} U_i$
are all simply connected.
Consider the groupoid
$${\mathcal H}=\coprod_{I\subseteq J} G_J\times U_J\rightrightarrows
\coprod_{I} G_I\times U_I
$$
where $I,J$ run over finite subsets of the index set. The structure maps are similar to those of the groupoid $\G$ constructed in section 3.
 The  groupoids ${\mathcal H}$ and $\G$  are easily seen to be equivalent in the sense of \cite[\S 2.4]{moerdijk}. Thus the weak homotopy type of
$B\mathcal{H}$ is also $[X]$.  Let $L$ be the locally 
constant sheaf of sets on $X$ corresponding to $\pi=\pi_1(X)$ with its left $\pi$-action.
We get a morphism of groupoids $\lambda:\mathcal{H}\to \pi$ with the elements $\lambda_{I,J}\in \pi$ corresponding to the  transition
functions of $L$ viewed as a flat bundle. The map $\lambda$ induces a a map of simplicial spaces $k_\dt: B\mathcal{H}_\dt\to B\pi_\dt$ whose geometric realization
is precisely $k$.

As noted in the introduction, theorem \ref{thm:splitting} is fairly elementary
when $X$ and $Y$ are smooth algebraic varieties. 

\begin{proof}[Proof  for algebraic varieties.]
Since the generic fibre $X_\eta$ has a rational point over some finite
extension of $\C(Y)$. There exists (by resolution of singularities)
a smooth  $Y'$ and a proper generically finite map $p:Y'\to Y$
such that $X\times_Y Y'\to Y'$ has a section. Therefore
$\pi_1(X\times_Y Y')\to \pi_1(Y')$ splits, and so the theorem
follows from lemma \ref{lemma:genfinred}.
\end{proof}

We employ a different strategy for the general case.

\begin{proof}[Proof of theorem \ref{thm:splitting}  in general]
By lemmas 3.1, 3.3, 3.4 and 3.5, we can reduce to the case
where $f:X\to Y$ is a connected holomorphic K\"ahler map of  orbifolds
satisfying the assumptions of lemma \ref{lemma:cko}.
Set $G= \pi_1(Y)$, $H= \pi_1(X)$ and $K =  \ker(\pi_1(f))$.
Let $V= H_1(K,\Q)= K/DK\otimes \Q$ with its natural $G$-action. We also
view this as a locally constant sheaf on $Y$. Then by
\cite[thm 4]{hs}, $e\otimes \Q$
is $\pm d_2(id)$, where
$$d_2:Hom_G(V,V)\cong H^0(G, H^1(K,V))\to H^2(G, H^0(K, V)) \cong H^2(G,V)$$
is the differential of the Hochschild-Serre spectral sequence.

After choosing compatible atlases for $X$ and $Y$, we can, for the purposes of sheaf
theoretic calculations, replace $f$
by a map of simplicial spaces $f_\dt:X_\dt\to Y_\dt$. More precisely, we can identify
sheaves on $X$ and $Y$ with certain simplicial sheaves on $X_\dt$ and $Y_\dt$.
In addition, the Leray spectral sequences for $f_\dt$ and  $f$  can be identified.
Consider the diagram
$$
\xymatrix{
 BH_\dt\ar[d]^{\phi_\dt} & X_\dt'\ar[d]^{F_\dt}\ar[l] & X_\dt\ar[d]^{f_\dt}\ar[l]\ar@/_/[ll]_{k_\dt} \\ 
 BG_\dt & Y_\dt\ar[l]^{k_\dt} & Y_\dt\ar[l]^{=}
}
$$
where the maps labeled by $k_\dt$ are the canonical maps realized simplicially as
above.
The left hand square is Cartesian. The geometric realization $\phi$ of $\phi_\dt$
can be assumed to be a fibration.  Consequently, the realization $F$ of $F_\dt$ is
also a fibration.

The Hochschild-Serre spectral sequence can
be identified with the Leray spectral for $\phi$ with coefficients
in the local system $\phi^*V$. This in turn can be identified with the
spectral sequence for $\phi_\dt$. 
Since the Leray spectral sequences
for $\phi_\dt,F_\dt$ and $f_\dt$, and hence $f$ are compatible, we have a commutative diagram
$$
\xymatrix{
 H^0(BG_\dt, R^1\phi_{\dt*}\phi_\dt^*V)\ar[r]\ar^{d_2}[d] & H^0(Y_\dt, R^1F_{\dt*}F_\dt^*V)\ar[r]\ar[d] & H^0(Y,R^1f_*f^*V)\ar^{d_2'}[d] \\ 
 H^2(BG_\dt,\phi_{\dt*}\phi_\dt^*V)\ar[r]^{k^*} & H^2(Y_\dt,F_{\dt*}F_\dt^*V)\ar[r]^{\ell} & H^2(Y,f_*f^*V)
}
$$
We can identify $k_\dt^*$ with the map on geometric realizations
$H^2(BG, \phi_*\phi^*V)\to H^2([Y],F_*F^*V)$.
This is injective, since the (homotopy) fibre of $k$ is simply connected.
Also since the fibres of the geometric realizations of $f_\dt$ and $F_\dt$
are connected, we have $F_{\dt*}F_\dt^*V=f_*f^*V=V$ by the
projection formula. Thus $\ell$ is an
isomorphism. So it suffices to prove that $d_2'$ is zero.

By our assumptions, we have an exact sequence
$$1\to L\to \pi_1(U)\to G\to 1$$
where the group $L$ is normal subgroup generated by powers of
loops $\gamma_i^{m_i}$
around components of the discriminant divisor $D$, and $U=Y-D$.
Since the first homology of the general fibre of $f$  surjects onto $V$ by lemma~\ref{lemma:cko},
the restriction
$V|_U$ can be identified with a locally constant quotient of $W=(R^1f_*\Q)^\vee|_U$ on which $L$ acts trivially.
$W$ carries a pure variation  of Hodge structure, therefore the action
of $\pi_1(U)$ on it is semisimple by Deligne's theorem \cite[7.2.5]{schmid}. Consequently $V|_U\subset
W$ is a direct summand. 
Therefore $V=IC(V)$ is a direct summand of
 $IC(W)=IC(f_*f^*W)$. Consequently, we
have a diagram
$$
\xymatrix{
 H^0(Y,R^1f_*f^*V)\ar^{d_2'}[r]\ar[d] & H^2(Y, f_*f^*V)\ar^{i}[d] \\ 
 H^0(Y, IC(R^1f_*f^*W))\ar^{d_2''}[r] & H^2(Y, IC(f_*f^*W))
}
$$
where the map labeled $i$ is injective. The map $d_2''$, which  is a summand
of the
differential of the perverse Leray spectral sequence, vanishes by
corollary \ref{cor:decomp}.
Therefore $d_2'=0$ and the theorem is proven.
\end{proof}

\begin{cor}\label{cor:discon}
  The theorem holds when $X$ is \K{} and $f$ is assumed to be proper surjective such that the
map $\pi_1(X)\to \pi_1(Y)$ is surjective. (Connectedness is not assumed.)
\end{cor}

\begin{proof}
By Stein factorization and resolution of singularities, we can find a commutative
diagram of complex manifolds
  $$
\xymatrix{
 X'\ar[r]\ar[d]^{f'} & X\ar[d]^{f} \\ 
 Y'\ar[r]^{g} & Y
}
$$
where $f'$ is as in the theorem and $g$ is generically finite. The result now follows
lemma~\ref{lemma:genfinred} and the theorem.
 \end{proof}

Let $\Gamma_g = \langle a_1,\ldots a_{2g}\mid [a_1,a_{g+1}]\ldots
[a_{g},a_{2g}]=1\rangle $
be the surface group of genus $g$. Call a homomorphism $h:\pi\to \Gamma_g$
maximal, if it does not factor through any $\Gamma_{g'}$ with $g'>g$.

\begin{cor}
  If a K\"ahler group admits a maximal surjective homomorphism $h:\pi\to \Gamma_g$ with $g\ge 2$, then $e(h)$ must be torsion.
\end{cor}

\begin{proof}
  By a theorem of Beauville-Siu \cite[thm 2.11]{abc}, $h$ is can be realized as
 $\pi_1(f)$ where $f:X\to C$ is a surjection onto a Riemann surface. This  uses
maximality. We can now finish the proof by appealing to the previous corollary, although in fact
it is not necessary. If $f$ were not  connected,
we could Stein factor the map, and conclude that $h$ is not maximal.
\end{proof}

This implies that a nontrivial central extension of $\Gamma_g$ by a torsion
free abelian group is not K\"ahler. In particular, this rules out the group
given in the introduction. J. Amor\'os has pointed that to me that one can also
see that this example has nontrivial Massey products, thus contradicting the formality
theorem of \cite{dgms}.

\begin{cor}
  Suppose that  $\pi$ is a K\"ahler group, and either $b_1(\pi)=\dim H^1(\pi,\Q)=2$,
or $b_1(\pi) =4$ and cup product  $\wedge^2 H^1(\pi,\Q)\to H^2(\pi,\Q)$ is injective.
Then $e(\alpha)$ is torsion, where $\alpha:\pi\to \pi/D\pi$ is the Abelianization.
\end{cor}

\begin{proof}
Let $X$ be a compact K\"ahler manifold with fundamental group $\pi$. 
The map $\pi\to \pi/D\pi/(\text{torsion})$ can be realized as the map 
as $\pi_1(alb)$, where $alb:X\to Alb(X)$ is the Albanese map.
We claim that $alb$ is surjective.
If $b_1(\pi)=2$, then $Alb(X)$ is an elliptic curve, and so surjectivity is clear.
When $b_1(\pi)=4$, $Alb(X)$ is a two dimensional torus.
So the image of $alb(X)$ is either a curve or all of $Alb(X)$.  The first
case is ruled out by the injectivity assumption for the cup product map. 
Therefore we are done by corollary~\ref{cor:discon}.
\end{proof}

\end{document}